\documentclass[review,onefignum,onetabnum]{Inputs/siamonline250211}

\newcommand{\Given}{\:|\:}

\newcommand{\mbX}{\mathbf{X}}

\newcommand{\mbf}{\mathbf{f}}

\newcommand{\mbg}{\mathbf{g}}

\newcommand{\mbx}{\mathbf{x}}
\newcommand{\mby}{\mathbf{y}}

\newcommand{\mbt}{\mathbf{t}}

\newcommand{\mbz}{\mathbf{z}}
\newcommand{\mbu}{\mathbf{u}}
\newcommand{\mbA}{\mathbf{A}}
\newcommand{\mbb}{\mathbf{b}}
\usepackage{mathrsfs}

\usepackage{lipsum}
\usepackage{amsfonts}
\usepackage{graphicx}
\usepackage{epstopdf}
\usepackage{algorithm}
\usepackage{algpseudocode}
\ifpdf
  \DeclareGraphicsExtensions{.eps,.pdf,.png,.jpg}
\else
  \DeclareGraphicsExtensions{.eps}
\fi

\usepackage{enumitem}
\setlist[enumerate]{leftmargin=.5in}
\setlist[itemize]{leftmargin=.5in}

\newsiamremark{remark}{Remark}
\newsiamremark{hypothesis}{Hypothesis}
\crefname{hypothesis}{Hypothesis}{Hypotheses}
\newsiamthm{claim}{Claim}
\newsiamremark{fact}{Fact}
\crefname{fact}{Fact}{Facts}

\headers{ANN Koopman eigenfunction solver}{M.~Kreider, P.~J.~Thomas, and Y.~Li}

\title{Artificial neural network solver for Fokker-Planck and Koopman eigenfunctions%\thanks{Submitted to the editors DATE.}
}

\author{%
  Max Kreider{\thanks{Department of Mathematics, Applied Mathematics, and Statistics,
            Case Western Reserve University.}
            \thanks{Corresponding author: \email{mbk62@case.edu}}}
  \and
  Peter J.~Thomas%
    \thanks{Department of Mathematics, Applied Mathematics, and Statistics,
            Case Western Reserve University.}% ← reuse “second” footnote
  \and
  Yao Li%
    \thanks{Department of Mathematics and Statistics, University of Massachusetts Amherst}%
}

\usepackage{amsopn}

\ifpdf
\hypersetup{
  pdftitle={Artificial neural network solver for Fokker-Planck and Koopman eigenfunctions},
  pdfauthor={M. Kreider, P.J. Thomas, and Y. Li}
}
\fi

\begin{document}

\maketitle

% REQUIRED
\begin{abstract}
For a stochastic differential equation (SDE) that is an It\^{o} diffusion or Langevin equation, the
Fokker-Planck operator governs the evolution of the probability density,  
while its adjoint, the infinitesimal generator of the stochastic Koopman operator, governs the evolution of system observables, in the mean.
The eigenfunctions of these operators provide a powerful framework to analyze SDEs, and have shown to be particularly useful for systems of stochastic oscillators.
However, computing these eigenfunctions typically requires solving high-dimensional PDEs on unbounded domains, which is numerically challenging.
Building on previous work, we propose a data-driven artificial neural network solver for Koopman and Fokker-Planck eigenfunctions.
Our approach incorporates the differential operator into the loss function, improving accuracy and reducing dependence on large amounts of accurate training data.
We demonstrate our approach on several numerical examples in two, three, and four dimensions.
\end{abstract}

% REQUIRED

% REQUIRED

\section{Introduction}

Many physical and biological systems have some degree of noise or uncertainty, and are widely modeled as systems of stochastic differential equations.
Traditionally, such systems have been characterized by studying the time-dependent and stationary probability densities associated with the underlying stochastic dynamics \cite{gardiner1985handbook,gu2021stationary,risken1996fokker}.
The Fokker-Planck operator describes the evolution of the probability density analytically.

However, recent attention has focused on the infinitesimal generator of the stochastic Koopman operator (SKO)\footnote{For brevity, we will write ``SKO'' to mean the \emph{infinitesimal generator} of the stochastic Koopman operator.} as a complementary descriptor of stochastic processes.
The SKO is adjoint to the Fokker-Planck operator, and  describes the evolution of the expectation of system observables.
Because the action of the SKO is linear on the space of $\mathcal{C}^2$-continuous system observables, the eigenspectrum of the SKO characterizes the original nonlinear, stochastic dynamics \cite{vcrnjaric2020koopman}.
SKO approaches have gained significant traction in fluid dynamics \cite{mezic2013analysis, rowley2009spectral}, control problems \cite{korda2018linear,otto2021koopman}, neuroscience \cite{brunton2016extracting,melland2023attractor}, and stochastic oscillators \cite{kreider2025q, Perez2023universal,perez2021isostables,thomas2014asymptotic}.

The space of system observables is infinite-dimensional, so obtaining finite dimensional approximations of the SKO is important for practical computations \cite{budivsic2012applied}.
While traditional PDE solvers (finite differences, finite elements) are effective for low-dimensional problems, they lose effectiveness for high-dimensional problems $(n\geq 3)$.
Non-trivial boundary conditions on unbounded domains can further complicate numerical implementation, as can the potential presence of continuous spectra \cite{mezic2005spectral}.
While much effort has been devoted to overcome these challenges, solving high-dimensional Fokker-Planck and Koopman eigenvalue problems is still a challenging task.

Machine learning is a powerful approach for solving high-dimensional PDEs and related eigenvalue problems.
Some works have employed diffusion Monte Carlo methods \cite{han2020solving, schatzle2021convergence}, deep Ritz methods \cite{ji2024deep,lu2021priori}, physics-informed neural network (PINN) \cite{raissi2019physics,sirignano2018dgm}, and BSDE methods \cite{han2020convergence, takahashi2022new, tsuchida2023control}, among other approaches \cite{barbara2023train,beck2020overview, heinlein2021combining,weinan2021algorithms}.
Recent work proposed a data-driven scheme that replaces boundary conditions by Monte Carlo simulation data when solving the Fokker-Planck equation \cite{dobson2019efficient,li2018data}.
This method was extended into an artificial neural network (ANN) approach, and was found to be effective for solving high-dimensional problems \cite{zhai2022deep}.
The idea behind the method is to train an ANN to approximate a probability density function that minimizes both (i) the distance to the Monte Carlo reference data, and (ii) the norm of the residual of the forward operator $\mathcal{L}$. Importantly, it was found that relatively inaccurate Monte Carlo approximations at only a small number of reference points were sufficient to guide the ANN to an accurate solution of the invariant distribution. We also refer to \cite{chen2021solving} for a PINN-based Fokker-Planck solver that uses a different regularization term based on Monte Carlo simulations. 

In this manuscript, we will extend previous results \cite{li2018data,zhai2022deep} to develop a novel ANN approach to solve for the low-lying eigenfunctions and eigenvalues of the SKO and Fokker-Planck operators associated with stochastic oscillators.
Leveraging results from \cite{thomas2014asymptotic}, we use Monte Carlo approximations of the time-dependent probability densities to estimate the low-lying eigenvalues of a given system.
Then, we obtain coarse approximations of the low-lying forward and Koopman eigenfunctions as solutions of a linear least squares problem.
This Monte Carlo data for the eigenmodes serves as a reference in a loss function (as a replacement for the boundary conditions, which are not explicitly considered), which is minimized with an ANN, analogous to the approach considered in \cite{zhai2022deep}.

Our ANN approach is advantageous for several reasons. 
Existing ANN solvers for eigenfunction mainly use normalization to avoid converging to the trivial solution \cite{han2020solving, jin2022physics}. 
In comparison, our approach directly estimates the value of an eigenfunction at numerous collocation points in the domain, leading to both fast convergence and improved solution quality. Other methods, such as extended dynamic mode decomposition (EDMD) and variants, have recently gained traction as data-driven methods to approximate the low-lying spectrum of the SKO \cite{brunton2016koopman,brunton2021modern,colbrook2023mpedmd, mezic2022numerical,korda2018convergence,williams2015data}.
However, such methods rely on choosing a dictionary of functions; a poor choice of dictionary results in poor performance of the EDMD algorithm, and dictionaries which are too large often result in computational bottlenecks and ill-conditioned systems.
In contrast, our novel ANN approach does not require \textit{a-priori} system knowledge, nor a collection of dictionary functions.
Instead, it requires only time-series data, making it more broadly applicable to a wide range of problems.
Our ANN approach is also capable of simultaneously approximating several low-lying eigenfunctions of the Fokker-Planck and Koopman operators.
Importantly, the user is free to define reference points in the domain at which solution values are required.
This mesh-free approach is advantageous for solving high-dimensional problems because it requires relatively few reference points.

The remainder of this work is organized in the following manner.
We first describe necessary assumptions and mathematical preliminaries in \S\ref{sec: ml : math prelim}.
Then, we detail our ANN approach and numerical algorithms in \S\ref{sec: ml algs}, discuss error associated with the eigenfunction approximation in \S\ref{sec: error decomp}, and apply our technique to several examples in two, three, and four dimensions in \S\ref{sec: ml examples}. 

\section{Mathematical preliminaries}
\label{sec: ml : math prelim}

We consider stochastic oscillators of the form
\begin{equation}\label{eq: ml : intro system}
    d\mathbf{X} = \mathbf{f}(\mathbf{X}) dt + \mathbf{g}(\mathbf{X})d\mathbf{W}(t)
\end{equation}
with $\mathbf{X}(t)$ a stochastic process taking values in a domain $\Omega\subseteq\mathbb{R}^n$, $\mbf:\Omega\to\mathbb{R}^n$ a $\mathcal{C}^2$-continuous vector field, $\mathbf{g}$ an $n\times m$ matrix with $\mathcal{C}^2$-continuous entries such that $\mbg\mbg^T$ is nonsingular for all $\mbx\in\Omega$, and $d\mathbf{W}(t)$ an $m \times 1$ vector of  increments of $m$ independent Wiener process. 
We interpret \eqref{eq: ml : intro system} in the sense of It\^{o}.

We assume throughout that the time-dependent transition density of the Markov process
\begin{equation}
    \rho(\mby,t\Given\mbx,s) = \lim_{|d\mby|\to 0}\frac{\mathbb{P}\left(\mbX(t)\in[\mby,\mby+d\mby)\Given \mbX(s)=\mbx\right)}{|d\mby|}
\end{equation}
exists for all $t>s$ and for all $\mbx$ and $\mby$ in $\Omega$.
The evolution of the  density associated with \eqref{eq: ml : intro system} satisfies the Fokker-Planck (forward) equation, which defines the forward operator $\mathcal{L}$:
\begin{equation}\label{eq: math prelims : forward_equation}
    \frac{\partial}{\partial t} \rho = \mathcal{L}[\rho] = -\sum_{i} \partial_{x_i}[\mathbf{f} \rho] + \frac{1}{2}\sum_{ij} \partial_{x_i x_j}\{ [\mathbf{g}\mathbf{g}^T]_{ij} \rho \}
\end{equation}
and also the backward equation, which defines the backward operator $\mathcal{L}^\dagger$:
\begin{equation}\label{eq: math prelims : backward_equation}
    -\frac{\partial}{\partial s} \rho = \mathcal{L}^\dagger[\rho] = \sum_{i} \mathbf{f}_i\partial_{y_i}[ \rho] + \frac{1}{2}\sum_{ij} [\mathbf{g}\mathbf{g}^T]_{ij}\partial_{y_i y_j}\{  \rho \}
\end{equation}
for $1\leq i,j \leq n$.
The operator $\mathcal{L}^\dagger$, aka the generator of the Markov process $\mbX$, is the infinitesimal generator of the stochastic Koopman operator \cite{vcrnjaric2020koopman}.
We henceforth refer to $\mathcal{L}^\dagger$ as ``the SKO'' for brevity.
We are concerned with the numerical solution of the forward and backward eigenvalue problems
\begin{equation}\label{eq: ml : eigenvalue equations}
    \begin{split}
        \mathcal{L}[P_\lambda(\mbx)] &= \lambda P_{\lambda}(\mbx)
        \\
        \mathcal{L}^\dagger [Q^*_\lambda(\mby)] &= \lambda Q^*_{\lambda}(\mby)
    \end{split}
\end{equation}
In particular, we focus on solving \eqref{eq: ml : eigenvalue equations} for the lowest-lying eigenpairs.

We assume that \eqref{eq: ml : intro system} is \textit{robustly oscillatory}: (i) all eigenvalues, apart from the simple trivial eigenvalue $\lambda_0\equiv 0$, have negative real parts, (ii) the first nontrivial eigenvalue, $\lambda_1 = \mu \pm i\omega $, is a unique complex-conjugate pair, with quality factor $|\omega/\mu| \gg 1$, and (iii) all other eigenvalues $\lambda'$ have more negative real parts, $\Re[{\lambda'}] < \upsilon\,\Re[{\lambda_1}]$, where $0<\upsilon\in \mathbb{R}$. 
We further assume that \eqref{eq: ml : intro system} has a discrete set of eigenvalues with biorthonormal eigenfunctions, and has a unique stationary distribution.
Under these assumptions, the forward transition density may be expressed as a linear combination of the forward and backward eigenmodes 
\begin{equation}\label{eq: ml : inf sum}
    \rho(\mathbf{y},t \Given \mathbf{x},s) = P_0 + \sum_{\lambda\neq 0} e^{\lambda(t-s)}P_\lambda(\mathbf{y})Q^*_\lambda(\mathbf{x})
\end{equation}
for times $t>s$.
In cases where there is a large spectral gap $\left|\Re[\lambda_1]-\sup\Re[\lambda']\right|$, previous works have truncated the infinite sum at only the first term \cite{ Perez2023universal, perez2021isostables, thomas2014asymptotic}
\begin{equation}\label{eq: ml : density expansion}
    \rho(\mathbf{y},t \Given \mathbf{x},s) - P_0(\mby) \approx    e^{\lambda_1(t-s)}P_{\lambda_1}(\mathbf{y})Q^*_{\lambda_1}(\mathbf{x}) + \text{c.c.}
\end{equation}
While our focus is on the $Q$-function, the unique complex lowest-lying eigenpair, in practice we achieve a more faithful approximation of system dynamics by including the first several eigenpairs, i.e., by truncating \eqref{eq: ml : inf sum} after several terms.
In the following, we will use \eqref{eq: ml : density expansion} to generate data-driven approximations of the forward and backward eigenfunctions.\footnote{Methods to extract approximations of the $Q$-function from time-series data were first introduced in \cite{thomas2014asymptotic} by Thomas and Lindner.} 
By assuming knowledge of high-order eigenvalues, i.e., including more terms using \eqref{eq: ml : inf sum}, we will obtain more accurate solutions for the $Q$-function.

Finally, we recall that the forward eigenfunctions should satisfy ``vanishing boundary conditions'' at infinity, and that the backward eigenfunctions should have vanishing derivatives at infinity, when $\Omega$ is unbounded.
In practice, it is not possible to implement these boundary conditions; instead, we restrict attention to a compact, rectangular computational domain, $\mathcal{R}\subset \Omega$.

\section{Numerical methods and ANN implementation}
\label{sec: ml algs}

We begin with an overview of our novel ANN approach to solve the eigenvalue problems \eqref{eq: ml : eigenvalue equations}.
First, realizations of the stochastic process \eqref{eq: ml : intro system} are used to generate Monte Carlo approximations of the time-dependent probability density.
The relaxation of the time-dependent probability density to an invariant distribution allows for an accurate estimation of the low-lying eigenvalue, $\lambda_1$, of the SKO and forward operator.
Then, we formulate a linear least squares problem using the eigenfunction expansion ansatz \eqref{eq: ml : density expansion}.
The solution of the least squares problem gives a coarse approximation of the forward and backward eigenmodes.
Finally, the least squares data is fed into an ANN, which interpolates and smooths the data and increases the accuracy of the approximation.

\subsection{Step 1: collocation points}

Here, the goal is to generate two sets of collocation points: $\mathcal{X}=\{\mbx_i, \; i=1,2,\dots,N_x\}$ and $\mathcal{Y} = \{\mby_j, \; j=1,2,\dots,N_y\}$.
We call $\mathcal{X}$  the ``training set'' and $\mathcal{Y}$  the ``reference set''.
Ultimately, a neural network will learn a specific eigenfunction on the training set to replace boundary conditions, and will learn the differential operator, $\mathcal{L}$ or $\mathcal{L}^\dagger$, on the reference set.
In other words, training on $\mathcal{Y}$ ensures that the output of the ANN satisfies (approximately) either the forward or backward equation, while training on $\mathcal{X}$ guides the ANN to a specific normalization of an eigenfunction and avoids the trivial solution.

To generate the training set, $\mathcal{X}$, we discretize the rectangular computational domain $\mathcal{R}$ into boxes of equal dimension, but we store in memory only $N_x$ of the boxes corresponding to the chosen collocation points.
Each box is represented by the coordinates corresponding to the center of the respective box.
This step is critical in higher dimensions, where storing a full mesh is not feasible. 
The points, $\mbx_i$, at the centers of the boxes will be associated with numerical estimates of the probability density.
In contrast, we do not bin the reference set in boxes because these points are used solely for minimizing the residual of the operator loss.

%evaluating the differential operators, and are not associated with numerical estimates for probability densities.

The support of the forward eigenmodes is often restricted to the invariant distribution, which may be concentrated near low dimensional manifolds or in specific regions of the domain. While to the best of our knowledge this is not supported by known results, we expect concentration results of invariant probability measures \cite{huang2018concentration, ji2019quantitative, li2016systematic} and quasi-stationary distributions \cite{shen2024concentration} to be applicable after some modifications.
In contrast, the support of the backward eigenmodes typically spans all of $\Omega$.
The generation of the collocation points should take this difference into account.

To generate collocation points, $\mathcal{X}$,
we first run \eqref{eq: ml : intro system} for a ``burn-in'' time, $T_{\text{burn}}$, so that the location of the terminal point of the process is effectively drawn from a distribution that is near the stationary distribution.
Then, we choose a ``ratio'', $\alpha\in [0,1]$, continue the realization of \eqref{eq: ml : intro system} starting from the terminal point of the burn-in simulation, and take only a fraction $\alpha$ (one should generally take $\alpha \in [0.5, 0.9]$) of the collocation points from this trajectory.
We establish a ``sample'' time, $t_{\text{gap}}$, and sample only at integer multiples of $t_{\text{gap}}$ so that points are not too close to each other.
The other $1-\alpha$ of the collocation points are chosen from a uniform distribution on the computational domain $\mathcal{R}$.
We associate each sampled point $\mbx_i$ with the box in $\mathcal{R}$ in which it lands provided it has not yet been taken by a collocation point. The process stops when $N_x$ collocation points are collected. We reiterate that only points $\mbx_i$ are binned in the boxes, and that one should take $\alpha=0$ when generating collocation points for the backward eigenmodes. 
This procedure ensures that many of the collocation points are in locations where the forward eigenmodes take on non-trivial values, and ultimately reduces the number of collocation points required for accurate estimation.

A similar procedure can be used to generate $\mathcal{Y}$, although these points are not binned in boxes because they are used to evaluate the differential operator and are not associated with numerical estimates of a probability density.
Algorithms SM1.1 and SM1.2 in the supplement provide more details for the generation of $\mathcal{Y}$ and $\mathcal{X}$, respectively \cite{zhai2022deep}.

\subsection{Step 2: density estimation}\label{step 2}

The goal of this step is to estimate time-dependent transition densities $\rho(\mathbf{y}, t \,|\, \mathbf{x}, s)$ associated with equation \eqref{eq: ml : intro system}.

To numerically approximate a time-dependent density, we fix a sampling time interval, $t_{\text{gap}}$, and the number of sample times, $N_t$.
Realizations of the process \eqref{eq: ml : intro system} are sampled at integer multiples of $t_{\text{gap}}$ so that samples are not too close together.\footnote{The parameter $t_{\text{gap}}$ should be  much smaller than the mean period of oscillation of a given system.}
These samples are recorded in a matrix $\mathcal{D}$, which has $N_t$ rows (one for each time slice) and $N_x$ columns (one for each collocation point, $\mbx_i\in\mathcal{X}$).

The procedures for estimating the densities associated with forward and backward eigenfunctions differ. To distinguish between them, we denote by $\mathcal{D}_F$ and $\mathcal{D}_B$ the matrices containing the estimated time-dependent densities for the forward and backward eigenfunctions, respectively. 

We begin by estimating the forward transition densities $\rho(\mbx_i, t_k \mid \mbx_0, 0)$ for a fixed initial condition $\mbx_0$, where $\mbx_i \in \mathcal{X}$ and $t_k \in {t_1, \dots, t_{N_t}}$. Using Monte Carlo simulation, we generate $K$ realizations of equation \eqref{eq: ml : intro system} starting from $\mbx_0$. At each sampling time $t_k$, we count how many trajectories fall into the bin centered at each collocation point $\mbx_i$. If the bin around $\mbx_i$ contains $K_i$ samples, the estimated transition probability is given by $K_i / (K \delta_i)$, where $\delta_i$ is the volume of the bin around $\mbx_i$. These estimated values are then stored in the matrix $\mathcal{D}_F$.

Letting $\tilde{\rho}_{f,i}(t_k) \equiv \tilde{\rho}(\mbx_i, t_k \mid \mbx_0, 0)$ denote the estimated transition density, it follows that the $i$th column of $\mathcal{D}_F$ provides an approximation to the time series ${\tilde{\rho}_{f,i}(t_k)}_{k=1}^{N_t}$. Since the initial condition $\mbx_0$ is fixed, the coefficient $Q_{\lambda_1}^*(\mbx_0)$ in the eigenfunction expansion ansatz \eqref{eq: ml : density expansion} is constant. If desired, one can also sample $\mbx_0$ from a prescribed distribution to reduce the influence of higher-order modes. Further implementation details can be found in Algorithm SM1.3 in the Supplement.

We now turn to the computation of the backward eigenmodes.
The objective in this step is to approximate the forward transition probability $\rho(\mathcal{B}, t \mid \mbx_i, 0)$ for each $\mbx_i \in \mathcal{X}$ at a sequence of time points $\mbt = \{t_k\}_{k=1}^{N_t}$.
Here, $\mathcal{B}$ denotes a reference box, $\mathcal{B} \subset \mathcal{R}$, which is chosen to contain a nontrivial portion of the invariant density.\footnote{The reference box is introduced as a practical approximation. Ideally, one would track the probability mass arriving in an individual box, but the likelihood of a trajectory landing in a single small bin is typically very low, necessitating a prohibitively large number of samples. To improve sampling efficiency, we instead consider a larger region $\mathcal{B}$.}
We denote the estimated transition probability into $\mathcal{B}$ at time $t_k$, given initial condition $\mbx_i$, by $\tilde{\rho}_{b,i}(t_k) \equiv \tilde{\rho}(\mathcal{B}, t_k \mid \mbx_i, 0)$.

To compute $\tilde{\rho}(\mathcal{B}, t_k \mid \mbx_i, 0)$, we simulate $K$ realizations of equation \eqref{eq: ml : intro system}, all initialized at $\mbx_i$. For each time $t_k$, we count the number of trajectories, say $K_j$, that lie within the set $\mathcal{B}$. The estimated transition probability is then given by $\tilde{\rho}_{b,i}(t_k) = K_j / K$, and the result is stored in the matrix $\mathcal{D}_B$. In our numerical examples, approximately $N_{\text{sample}} = 10^4$ trajectories per $\mbx_i$ are sufficient to achieve accurate estimates. Additional implementation details are provided in Algorithm SM1.4 of the Supplement. It is important to note that in estimating $\tilde{\rho}_{b,i}(t_k)$, the terminal region $\mathcal{B}$ is fixed. As a result, in the eigenfunction expansion ansatz \eqref{eq: ml : density expansion}, the integrals $\int_{\mathcal{B}} d\mby \; P_{\lambda_1}(\mby)$ and $\int_{\mathcal{B}} d\mby\; P_{0}(\mby)$ remain constant across all initial conditions $\mbx_i$. This property is exploited in Step 4 (see \S\ref{ml: step 4}) to recover the backward eigenfunctions over the domain $\mathcal{X}$.

For low-dimensional problems, standard grid-based methods may be used (e.g., see \cite{li2018data}) to sample the probability density at the points $\mbx_i$.
For higher dimensional problems, a grid-based approach is too slow and memory intensive.
Instead, we opt for the memory-efficient Monte Carlo sampling algorithm introduced in \cite{zhai2022deep} (see Appendix B in that work for a complete description).

\subsection{Step 3: Eigenvalue Approximation}\label{step 3}

In this section, the goal is to estimate the low-lying eigenvalue, $\lambda_1=\mu_1 + i\omega_1$, corresponding to the $Q$-function and related forward eigenfunction.
The eigenvalue is then used to approximate the corresponding eigenfunction in Step \ref{ml: step 4}.

The eigenvalue approximation relies on one of the numerically estimated densities from Step 2.
Recall that the estimated densities arising from Algorithm SM1.3 are denoted by $\tilde{\rho}_{f,i}(t_k)$, and the estimated densities arising from Algorithm SM1.4 are denoted by $\tilde{\rho}_{b,i}(t_k)$.
As before, $\mbt=\{t_k\}_{k=1}^{N_t}$ is understood to represent the $N_t$-dimensional vector of the time slices.
For concreteness, we present results in terms of the $i$th density, $\tilde{\rho}_{f,i}(t_k)$.\footnote{One should choose $i$ such that $\mbx_i$ is near a nontrivial portion of the invariant density.
It is also practical to choose a reference box centered on a region with substantial  invariant density and sum together each $\tilde{\rho}_{f,i}(t_k)$ corresponding to locations in the reference box to reduce noise.
The exact choice of box will influence $b_1$, $b_3$, and $b_5$ in \eqref{eq: ml : objective function eval}, but will not influence the estimation of $\lambda_1$ (in the limit of large data), because $\lambda_1$ is a global property of the system.}

Because the underlying system is a stochastic oscillator, we expect $\tilde{\rho}_{f,i}(t_k)$ to have the form of a decaying oscillation.
Our goal is to estimate the decay rate (real part of $\lambda_1$) and frequency (imaginary part of $\lambda_1$) of this oscillation.
To that end, we define a function
\begin{equation}\label{eq: ml : objective function eval}
    f(t;\mathbf{b}) = b_1\exp(b_5 t)\sin\left( \frac{2\pi t}{b_2} + \frac{2\pi}{b_3}\right) + b_4
\end{equation}
with  parameters $\mathbf{b}=[b_1,b_2,b_3,b_4,b_5]$ to be determined.
This function is chosen to model both contraction and oscillatory behavior, where:
\begin{itemize}
    \item $b_1$ and $b_4$ adjust the amplitude and vertical shift,
    \item $b_5$ captures system damping, and approximates the real part of the eigenvalue,
    \item $b_2$ and $b_3$ describe the periodicity and phase shift of the oscillation, where $2\pi/b_2$ approximates the imaginary part of the eigenvalue.
    \item $b_1$, $b_2$, $b_4$ and $-b_5$ are strictly positive.
\end{itemize}
Then, we define an objective function
\begin{equation}\label{eq: eigenvalue lsqr}
    F(\mbt;\mathbf{b}) = \sum_{j=t_s}^{t_f} (\tilde{\rho}_{f,i}(t_j) - f(t_j;\mathbf{b}))^2
\end{equation}
which measures the least squares error between the data, $\tilde{\rho}_{f,i}(t_k)$, and the model.
By solving \eqref{eq: eigenvalue lsqr} in the (nonlinear) least squares sense\footnote{We use MATLAB's ``fminsearch'' command to perform the fitting.}, we obtain estimates for the parameters $\mathbf{b}$.
The parameter $b_5$ approximates the real part of $\lambda_1$, while $2\pi/b_2$ approximates the imaginary part of $\lambda_1$.
We perform this fitting for different collocation points, and let $\tilde{b}_5$ and $\tilde{b_2}$ represent the average output of the fitting procedure.
We let $\Tilde{\lambda}_1=\Tilde{\mu}_1+i\Tilde{\omega}_1 = \tilde{b}_5 + (2\pi/\tilde{b}_2)i$ denote the least squares estimate of the true eigenvalue, $\lambda_1$.

Notice that the sum is not necessarily taken over all the sample points $t_k$.
Instead, it is taken over an intermediate time interval $t_1 \leq t_s <  t_f \leq t_{N_t}$. The formulation \eqref{eq: eigenvalue lsqr} assumes that significant contributions to the decaying oscillation arise only from the slowest decaying eigenvalue, $\lambda_1$.
However, in practice, even if the spectral gap (recall the \textit{robustly oscillatory criteria} in \S\ref{sec: ml : math prelim}) is large, we find that the effects of higher-order eigenmodes cannot be ignored. 
As will be discussed in \S\ref{sec: error decomp}, large $t_s$ reduces the effect of higher-order eigenmodes but increases the Monte Carlo error. Therefore, the intermediate time interval $[t_s,t_f]$ should be chosen so that the contributions from the higher-order modes decay but the decayed oscillation has not been overwhelmed by noise. 

%Therefore, the intermediate time interval $[t_s,t_f]$ should be chosen so that the contributions from the higher-order modes decay, but also so that the contribution from the slowest decaying mode has not relaxed to the stationary distribution.
%In other words, we expect the ansatz \eqref{eq: ml : density expansion} to be valid only for intermediate time intervals.

In some cases, it is useful to obtain estimates for the eigenvalues of higher-order eigenmodes.
If the stochastic oscillator reduces to a stable limit cycle oscillator in the vanishing noise limit, the low-lying eigenvalue spectrum may be approximately  parabolic with $\lambda_n\approx \mu_1 n^2 + i\omega_1 n$ for integer $n\in \mathbb{Z}$.
This structure is consistent with the existence of a change of coordinates under which system dynamics may be approximated as a diffusion on the ring: $d\phi=\omega_1 dt + \sqrt{-2\mu_1}\,dW(t)$, whose eigenvalues lie exactly on this parabola
\cite{Houzelstein2024generalized,thomas2014asymptotic}.
In such cases the approximated eigenvalue $\Tilde{\lambda}_1$ may be used to estimate eigenvalues corresponding to higher-order modes.

\subsection{Step 4: Eigenfunction approximation}\label{ml: step 4}

We now discuss the approximation of the $Q$-function and the corresponding forward eigenfunction.
In what follows, we assume knowledge of the stationary distribution, $P_0$, as the stationary distribution can be accurately computed by using methods in \cite{dobson2019efficient} or \cite{zhai2022deep}.
We also require the eigenfunction expansion ansatz \eqref{eq: ml : density expansion},  which we reproduce here for convenience (setting $s=0<t$)
\begin{equation}\label{eq: alg : density expansion}
    \rho(\mathbf{y},t \Given \mathbf{x},0) - P_0(\mby) \approx    e^{\lambda_1 t}P_{\lambda_1}(\mathbf{y})Q^*_{\lambda_1}(\mathbf{x}) + \text{c.c.} = 2\Re\left[e^{\lambda_1 t}P_{\lambda_1}(\mathbf{y})Q^*_{\lambda_1}(\mathbf{x})\right]
\end{equation}

\paragraph{Forward eigenfunction}
To numerically approximate the forward eigenfunction, we consider the ansatz
\begin{equation}\label{eq: ml : forward approx 1}
    \tilde{\rho}_{f,i}(t_j) - P_0(\mbx_i) = 2\Re[\exp(\tilde{\lambda}_1 t_j) P_{\lambda_1}(\mbx_i) Q_{\lambda_1}^*(\mbx_0)]
\end{equation}
where we replace the continuous terminal variable $\mby$ with the collocation points: $\mby \to \mbx_i \in \mathcal{X}$, replace the continuous initial variable $\mbx$ with points sampled from the starting box: $\mbx \to \mbx_0$, and replace continuous time with discrete time slices: $t\to t_j$.
Here, $\{t_j\}_{j=T_s}^{T_f}\subset \{t_k\}_{k=1}^{N_t}$  is a subset of all the time-slices with $t_1 \leq t_{T_s} < t_{T_f} \leq t_{N_t}$.
We also replace the unknown true transition density with the known numerically approximated density, $\tilde{\rho}_{f,i}(t_j)$, and the unknown true eigenvalue with the approximated eigenvalue $\Tilde{\lambda}_1$.
Notice that, as in Step 3 concerning the approximation of the eigenvalue $\lambda_1$, we consider only a subset of the possible time-slices due to contamination from higher-order modes.

Note that in \eqref{eq: ml : forward approx 1}, only the quantities $P_{\lambda_1}(\mbx_i)$ and $Q_{\lambda_1}^*(\mbx_0)$ are unknown.
Our goal is to approximate $P_{\lambda_1}(\mbx_i)$ for each $\mbx_i\in \mathcal{X}$.
To that end, recall that $\tilde{\rho}_{f,i}(t_k)\equiv \tilde{\rho}(\mbx_i, t_k \Given \mbx_0, 0)$ where $\mbx_0$ is a fixed initial condition.
Consequently, the quantity $Q_{\lambda_1}^*(\mbx_0)$ remains constant across all indices $i$, and does not depend on the terminal point $\mbx_i$.
We absorb this constant factor into the definition of a rescaled eigenfunction
\begin{equation}
    \tilde{P}_{\lambda_1}(\mbx_i) = 2Q_{\lambda_1}^*(\mbx_0) P_{\lambda_1}(\mbx_i)
\end{equation}
With this notation, we rewrite \eqref{eq: ml : forward approx 1} as
\begin{equation}\label{eq: ml : forward approx 2}
    \tilde{\rho}_{f,i}(t_j) - P_0(\mbx_i) = \Re[\exp(\tilde{\lambda}_1 t_j) \tilde{P}_{\lambda_1}(\mbx_i)]
\end{equation}
We find it convenient to introduce the notation $\tilde{P}_{\lambda_1} = P^R_{\lambda_1}+i P^I_{\lambda_1}$ to express \eqref{eq: ml : forward approx 2} as
\begin{equation}\label{eq: ml : forward equation}
\begin{split}
    \tilde{\rho}_{f,i}(t_j) -  P_0(\mbx_i) - e^{\Tilde{\mu}_1 t_j}\Bigg[&\cos(\Tilde{\omega}_1 t_j)P_{\lambda_1}^R(\mbx_i) - \sin(\Tilde{\omega}_1 t_j)P_{\lambda_1}^I(\mbx_i)\Bigg] = 0 
\end{split}
\end{equation}
which is a linear system $\mbA_f \mbu_f = \mbb_f$, with
\begin{equation}
   \mbb_f =  \begin{bmatrix}
        \exp(-\tilde{\mu}_1 t_{T_s}) [\tilde{\rho}_{f,i}(t_{T_s}) -  P_0(\mbx_i)]
        \\
        \exp(-\tilde{\mu}_1 t_{T_s+1})[\tilde{\rho}_{f,i}(t_{T_s+1}) -  P_0(\mbx_i)]
        \\
        \vdots
        \\
        \exp(-\tilde{\mu}_1 t_{T_f})[\tilde{\rho}_{f,i}(t_{T_f}) -  P_0(\mbx_i)]
    \end{bmatrix}, \quad \mbu_f = \begin{bmatrix}
        P_{\lambda_1}^R(\mbx_i)
        \\
        P_{\lambda_1}^I(\mbx_i)
    \end{bmatrix}
\end{equation}
and
\begin{equation}
    \mbA_f = \begin{bmatrix}
        \cos(\tilde{\omega}_1 t_{T_s}) &  -\sin(\tilde{\omega}_1 t_{T_s})
        \\
        \cos(\tilde{\omega}_1 t_{T_s+1}) &  -\sin(\tilde{\omega}_1 t_{T_s+1})
        \\
        \vdots 
        \\
        \cos(\tilde{\omega}_1 t_{T_f})&  -\sin(\tilde{\omega}_1 t_{T_f})
    \end{bmatrix}
\end{equation}
Equation \eqref{eq: ml : forward equation} is a linear least squares problem, and should be solved separately for each $\mbx_i\in \mathcal{X}$.\footnote{We use MATLAB's ``backslash'' to perform this computation.}
The solution of \eqref{eq: ml : forward equation} gives an approximation of the rescaled forward eigenfunction, $\tilde{P}_{\lambda_1}(\mbx_i)$ for  $\mbx_i \in \mathcal{X}$.

We note that the rescaled function \( \tilde{P}_{\lambda_1} \) does not, in general, satisfy a biorthogonality condition as it absorbs the unknown scalar factor \( 2Q_{\lambda_1}^*(\mbx_0) \).\footnote{If $Q_{\lambda_1}^*(\mbx_0)=0$, then our estimation procedure fails.
However, the $Q$-function for a robustly oscillatory stochastic oscillator typically has vanishing magnitude only at isolated points, namely at ``phaseless points" that tend to lie close to deterministic phaseless points (fixed points of the deterministic system).
One approach to ensure that $Q_{\lambda_1}^*(\mbx_0)\neq0$ is to choose $\mbx_0$ near a high density region of the invariant distribution.}
As a result, \( \tilde{P}_{\lambda_1} \) approximates the true forward eigenfunction \( P_{\lambda_1} \) only up to a complex constant.
If normalization is required, e.g., to enforce biorthogonality, this scalar can be estimated by computing the inner product \( \langle \tilde{P}_{\lambda_1}, Q_{\lambda_1}^* \rangle \), provided that an approximation of \( Q_{\lambda_1}^* \) is available. 
The function \( \tilde{P}_{\lambda_1} \) can then be rescaled accordingly to recover a biorthonormal representative of the forward eigenfunction.

\paragraph{Backward eigenfunction} The computational procedure to approximate the backward eigenfunction is similar to that of the forward case.
We consider the ansatz
\begin{equation}\label{eq: ml : backward approx 1}
    \tilde{\rho}_{b,i}(t_j) - \int_{\mathcal{B}} d\mby \; P_0(\mby) = 2\Re\left[\exp(\tilde{\lambda}_1 t_j) \left(\int_{\mathcal{B}} d\mby \; P_{\lambda_1}(\mby) \right) Q_{\lambda_1}^*(\mbx_i)\right]
\end{equation}
where we integrate the continuous terminal variable $\mby$ over the reference box $\mathcal{B}$, replace the continuous initial variable $\mbx$ with the starting boxes: $\mbx \to \mbx_i\in \mathcal{X}$, and replace continuous time with discrete time-slices: $t\to t_j$ with $t_j$ as above.
We also replace the true unknown density with the known numerically approximated density, $\tilde{\rho}_{b,i}(t_j)$, and the unknown true eigenvalue with the approximated eigenvalue $\Tilde{\lambda}_1$.

Note that in \eqref{eq: ml : backward approx 1}, only the quantities $\int_{\mathcal{B}} d\mby \; P_{\lambda_1}(\mby)$ and $Q_{\lambda_1}^*(\mbx_i)$ are unknown.
Our goal is to approximate $Q_{\lambda_1}^*(\mbx_i)$ for each $\mbx_i\in \mathcal{X}$.
To that end, recall that $\tilde{\rho}_{b,i}(t_k)\equiv \tilde{\rho}(\mathcal{B}, t_k \Given \mbx_i, 0)$ where $\mathcal{B}$ is a fixed reference box for each initial condition $\mbx_i$, and is used as an approximation of a terminal box near the invariant distribution.
Consequently, the quantity $\int_{\mathcal{B}} d\mby \; P_{\lambda_1}(\mby)$ remains constant across all indices $i$, and does not depend on the starting point $\mbx_i$.
We absorb this constant factor into the definition of a rescaled eigenfunction
\begin{equation}
    \tilde{Q}_{\lambda_1}^*(\mbx_i) = 2\left(\int_{\mathcal{B}} d\mbz \; P_{\lambda_1}(\mbz)\right) Q_{\lambda_1}^*(\mbx_i)
\end{equation}
With this notation, we rewrite \eqref{eq: ml : backward approx 1} as
\begin{equation}\label{eq: ml : backward approx 2}
    \tilde{\rho}_{b,i}(t_j) - \int_{\mathcal{B}} d\mbz \; P_0(\mbz) = \Re[\exp(\tilde{\lambda}_1 t_j) \tilde{Q}_{\lambda_1}^*(\mbx_i)]
\end{equation}
We find it convenient to introduce the notation $\tilde{Q}_{\lambda_1}^* = Q^R_{\lambda_1}+i Q^I_{\lambda_1}$ to express \eqref{eq: ml : backward approx 2} as
\begin{equation}\label{eq: ml : backward equation}
\begin{split}
    \tilde{\rho}_{b,i}(t_j) -  \int_{\mathcal{B}} d\mbz \; P_0(\mbz) - e^{\Tilde{\mu}_1 t_j}\Bigg[&\cos(\Tilde{\omega}_1 t_j)Q_{\lambda_1}^R(\mbx_i) - \sin(\Tilde{\omega}_1 t_j)Q_{\lambda_1}^I(\mbx_i)\Bigg] = 0 
\end{split}
\end{equation}
which is a linear system $\mbA_b \mbu_b = \mbb_b$, with
\begin{equation}
   \mbb_b =  \begin{bmatrix}
        \exp(-\tilde{\mu}_1 t_{T_s}) [\tilde{\rho}_{b,i}(t_{T_s}) -  \int_{\mathcal{B}} d\mbz \; P_0(\mbz)]
        \\
        \exp(-\tilde{\mu}_1 t_{T_s+1})[\tilde{\rho}_{b,i}(t_{T_s+1}) -  \int_{\mathcal{B}} d\mbz \; P_0(\mbz)]
        \\
        \vdots
        \\
        \exp(-\tilde{\mu}_1 t_{T_f})[\tilde{\rho}_{b,i}(t_{T_f}) -  \int_{\mathcal{B}} d\mbz \; P_0(\mbz)]
    \end{bmatrix}, \quad \mbu_b = \begin{bmatrix}
        Q_{\lambda_1}^R(\mbx_i)
        \\
        Q_{\lambda_1}^I(\mbx_i)
    \end{bmatrix}
\end{equation}
and
\begin{equation}
    \mbA_b = \begin{bmatrix}
        \cos(\tilde{\omega}_1 t_{T_s}) &  -\sin(\tilde{\omega}_1 t_{T_s})
        \\
        \cos(\tilde{\omega}_1 t_{T_s+1}) &  -\sin(\tilde{\omega}_1 t_{T_s+1})
        \\
        \vdots 
        \\
        \cos(\tilde{\omega}_1 t_{T_f})&  -\sin(\tilde{\omega}_1 t_{T_f})
    \end{bmatrix}
\end{equation}
We use different notation for $\mbA_f$ and $\mbA_b$ to highlight that different intermediate time intervals may be used for the approximation of the forward and backward eigenfunctions, respectively. 
Equation \eqref{eq: ml : backward equation} is a linear least squares problem, and should be solved separately for each $\mbx_i\in \mathcal{X}$.
The solution of \eqref{eq: ml : backward equation} gives an approximation of the rescaled backward eigenfunction, $\tilde{Q}_{\lambda_1}^*(\mbx_i)$ for  $\mbx_i \in \mathcal{X}$.

We note that, as in the case of the forward eigenfunction, the rescaled function \( \tilde{Q}_{\lambda_1}^* \) does not, in general, satisfy a biorthogonality condition as it absorbs the unknown scalar factor \( 2\int_{\mathcal{B}} d\mbz \; P_{\lambda_1}(\mbz) \).\footnote{Note that this constant factor is non-vanishing provided that the reference box overlaps  a non-vanishing region of the invariant distribution.}
As a result, \( \tilde{Q}_{\lambda_1}^* \) approximates the true backward eigenfunction \( Q_{\lambda_1}^* \) only up to a complex constant.
If needed, appropriate normalization can be enforced in a manner similar to that described for the case of the forward eigenfunction.

\paragraph{Computation of higher order eigenmodes} In some cases where the spectral gap is small, including more than one eigenmode in the least squares formulation can improve accuracy and allow for the approximation of higher-order eigenmodes, even if only the first eigenfunction is used in the artificial neural network training. 
We refer  \S\ref{sec: error decomp} for more analysis on how the linear regression error depends on higher-order modes.
To account for higher-order eigenmodes, we expand \eqref{eq: alg : density expansion} to include additional terms.
We now list the augmented forward and backward least squares problems that incorporate $M$ eigenmodes:
\begin{equation}
\label{eq-3-16}
\begin{split}
    \tilde{\rho}_{f,i}(t_j) -  P_0(\mbx_i) &- \sum_{k=1}^M e^{\Tilde{\mu}_k t_j}\Bigg[\cos(\Tilde{\omega}_k t_j)P_{\lambda_k}^R(\mbx_i) - \sin(\Tilde{\omega}_k t_j)P_{\lambda_k}^I(\mbx_i)\Bigg] = 0 
\end{split}
\end{equation}
and 
\begin{equation}
\label{eq-3-17}
\begin{split}
    \tilde{\rho}_{b,i}(t_j) -  \int_{\mathcal{B}} d\mbz \; P_0(\mbz) &- \sum_{k=1}^M e^{\Tilde{\mu}_k t_j}\Bigg[\cos(\Tilde{\omega}_k t_j)Q_{\lambda_k}^R(\mbx_i) - \sin(\Tilde{\omega}_k t_j)Q_{\lambda_k}^I(\mbx_i)\Bigg] = 0 
\end{split}
\end{equation}
respectively.
We numerically investigate the impact of including higher order modes in \S\ref{sec: ml examples}.
Note that our theory holds if a higher-order eigenmode and corresponding eigenvalue are purely real, as the imaginary parts are simply zero in the above expressions. %\yl{why our theory only holds if all eigenvalues are real? I thought we assume real eigenfunctions but eigenvalues are complex.}

\subsection{Step 5: ANN implementation}\label{ml: step 5}

The approximations of the eigenfunctions from step 4 are only available at collocation points $\mbx_i\in \mathcal{X}$.
Moreover, because these approximations rely on Monte Carlo simulations of the underlying probability density, the approximations are often noisy, especially in high dimensions.
Here, our goal is to use an ANN to generate eigenfunction approximations at any point of interest in $\mathcal{R}$, and to reduce noise.

Building on the approach pioneered by \cite{zhai2022deep}, we consider an objective function corresponding to operator $\mathscr{L}$ with eigenpair $(\lambda_n,U_{\lambda_n})$:
\begin{equation}
\label{eq: ml : object}
    \mathcal{F}_{\mathscr{L}}(\lambda_n,U_{\lambda_n}) =  \|(\mathscr{L}-\lambda_n)U_{\lambda_n}\|^2_{L^2(\mathcal{R})}  
    + \|U_{\lambda_n}-\Tilde{U}_{\lambda_n}\|^2_{L^2(\mathcal{R})}
\end{equation}
Here, the generic operator $\mathscr{L}$ may play the role of either the forward operator, $\mathcal{L}$, or the backward operator, $\mathcal{L}^\dagger$.
Accordingly, the eigenfunction $U_{\lambda_n}$ represents either a forward eigenfunction, $P_{\lambda_n}$, or a backward eigenfunction, $Q^*_{\lambda_n}$.
The notation $\Tilde{U}_{\lambda_n}$ is understood to represent a least-squares approximation of either the forward or backward eigenfunction from Step 4.

The first term in \eqref{eq: ml : object} ensures that $U_{\lambda_n}$ is a solution of the differential operator, i.e., one of the infinitely many eigenfunctions which are solutions to the eigenvalue problem without boundary conditions.
The last term circumvents the need for boundary conditions by enforcing a particular normalization via  reference data.

Our goal is to use an ANN to minimize the first term of the objective function \eqref{eq: ml : object} over the reference points, $\mathcal{Y}$, and to minimize the last term of the objective function over the training set, $\mathcal{X}$.
To that end, consider an ANN with parameter set $\chi$ (further details are provided in the supplement) with two outputs, $\mathcal{N}^R(\mbx;\chi)$ and $\mathcal{N}^I(\mbx;\chi)$, corresponding to the real and imaginary parts of an eigenfunction, respectively.
To approximate the eigenpair $(\lambda_n,U_{\lambda_n})$, we introduce the loss function $J_{\mathscr{L}}(\chi) \equiv J_{\mathscr{L}}(\tilde{\lambda}_n,\tilde{U}_{\lambda_n}(\mathcal{X}),\mathcal{Y};\chi)$:\footnote{We break the loss function into real and imaginary parts because TensorFlow 2 cannot work directly with complex numbers.}
\begin{equation}
    \begin{split}
        J_{\mathscr{L}}(\chi) &=  \frac{1}{N_y} \sum_{j=1}^{N_y} \left(\mathscr{L}[\mathcal{N}^R(\mby_j;\chi)] - \Tilde{\mu}_n\mathcal{N}^R(\mby_j;\chi) + \Tilde{\omega}_n\mathcal{N}^I(\mby_j;\chi) \right)^2
    \\
    &\quad+ \frac{1}{N_y} \sum_{j=1}^{N_y} \left(\mathscr{L}[\mathcal{N}^I(\mby_j;\chi)] - \Tilde{\mu}_n\mathcal{N}^I(\mby_j;\chi) - \Tilde{\omega}_n\mathcal{N}^R(\mby_j;\chi) \right)^2
    \\
    &\quad+ \frac{1}{N_x}\sum_{i=1}^{N_x} \left( \mathcal{N}^R(\mby_j;\chi) - \Re[\Tilde{U}_{\lambda_n}(\mbx_i)] \right)^2 + \frac{1}{N_x}\sum_{i=1}^{N_x} \left( \mathcal{N}^I(\mby_j;\chi) - \Im[\Tilde{U}_{\lambda_n}(\mbx_i)] \right)^2
    \end{split}
\end{equation}
Minimization of the loss function $J_{\mathscr{L}}(\chi)$ with respect to $\chi$ generates parameters, $\chi^*$, for which the network solution approximately satisfies 
\begin{equation}
    \mathscr{L}[\mathcal{N}^R(\mbx;\chi^*)+i\mathcal{N}^I(\mbx;\chi^*)] = \lambda_n [\mathcal{N}^R(\mbx;\chi^*)+i\mathcal{N}^I(\mbx;\chi^*)]
\end{equation}
As noted in \cite{zhai2022deep}, our setup is similar to the so-called Physics-Informed Neural Network (PINN) \cite{cuomo2022scientific}, in which the first term  over $\mathcal{Y}$ in the loss function (broken into real and imaginary parts) learns the ``physics'' of the problem via the differential operator, and the last term over $\mathcal{X}$ (broken into real and imaginary parts) is similar to boundary conditions.

The neural network training follows the ``Alternating Adam" approach used in \cite{zhai2022deep}. In each iteration, we shuffle two small batches from $\mathcal{X}$ and $\mathcal{Y}$ respectively, and update the parameters $\chi$ using the Adam optimizer. This is because the loss terms with the differential operators could have very different magnitudes than the loss terms corresponding to the reference data. The Adam optimizer is also almost invariant to rescaling \cite{kingma2014adam}.

%\pt{(*** PT stopped here 5/28/2025 ***)}

\section{Error decomposition of eigenfunction approximation}\label{sec: error decomp}

In this section we will study the error decomposition when estimating the forward eigenfunctions $P_{\lambda}(\mathbf{x})$ by solving a least squares problem, where $\mathbf{x}$ is a given collocation point. 
The case of the backward eigenfunctions is analogous. 
We start from equations \eqref{eq-3-16} and \eqref{eq-3-17} by assuming that the first $M$ eigenvalues denoted by $\lambda_1 = \mu_1 + i\,\omega_1 , \cdots, \lambda_M = \mu_M + i \,\omega_N $ (along with their complex conjugates), are known. 
This assumption helps us to focus on the effect of high-order modes, which is the main source of the error term in the eigenfunction. 
The first $M$ known eigenvalues can either come from the nonlinear regression method introduced in \S\ref{step 3} or other approaches like the coupling method \cite{dobson2021using}, EDMD \cite{brunton2016koopman,brunton2021modern,colbrook2023mpedmd, korda2018convergence, mezic2022numerical,williams2015data}, or the Matrix Pencil Method \cite{hua1990matrix, laroche1993use, sarkar1995using}.

Let $\{t_j\}_{j = 1}^N$ be the sampling times with $t_1 \equiv  t_{T_s}$ and $t_N \equiv t_{T_f}$ as in \S\ref{ml: step 4}. 
With the first $M$ eigenvalues known, the eigenfunction estimation problem becomes a linear least square problem for each $\mbx_i\in \mathcal{X}$
\begin{equation}
    \label{eqn:LSQ}
    \min_{\beta \in \mathbb{R}^{2M}} \| \mathscr{X} \beta - \mathscr{Y}\|^2
\end{equation}
for 
\begin{equation}
    \begin{split}
        \mathscr{X} &= \begin{bmatrix}
    e^{\mu_1 t_1}\cos (\omega_1 t_1) & e^{\mu_1 t_1}\sin (\omega_1 t_1) &\cdots & \cdots& e^{\mu_M t_1}\cos (\omega_M t_1) & e^{\mu_M t_1}\sin ( \omega_M t_1)
    \\
    e^{\mu_1 t_2}\cos ( \omega_1 t_2) & e^{\mu_1 t_2}\sin( \omega_1 t_2)&\cdots & \cdots & e^{\mu_M t_2}\cos( \omega_M t_2)& e^{\mu_M t_2}\sin( \omega_M t_2)
    \\
    \vdots &\vdots &\vdots &\vdots &\vdots &\vdots
    \\
    e^{\mu_1 t_N}\cos( \omega_1 t_N)& e^{\mu_1 t_N}\sin( \omega_1 t_N)&\cdots & \cdots &e^{\mu_M t_N}\cos( \omega_M t_N)& e^{\mu_M t_N}\sin (\omega_M t_N)
\end{bmatrix}
\\
\mathscr{Y} &= [\tilde{\rho}_{f,i}(t_1) - P_0(\mathbf{x}), \cdots, \tilde{\rho}_{f,i}(t_N) - P_0(\mathbf{x})]^T 
    \end{split}
\end{equation}
Define $\beta_0$ as the estimated value of the ground truth
\begin{equation}
    \beta_0 := [P^R_{\lambda_1}(\mathbf{x_i}), P^I_{\lambda_1}(\mathbf{x_i}), \cdots, P^R_{\lambda_M}(\mathbf{x_i}), P^I_{\lambda_M}(\mathbf{x_i})]^T 
\end{equation}
Assume further that $\mathscr{Y}$ has a decomposition
\begin{equation}
    \label{y-decomposition}
    \mathscr{Y} = \mathscr{X}^T\beta_0 + \mathscr{Y}_m + \mathscr{Y}_h
\end{equation}
where $\mathscr{Y}_h$ is the high-mode term satisfying
\begin{equation}
    \mathscr{Y}_h = C_1[e^{\hat\mu t_1}\cos (\hat \omega t_1 + \phi), \cdots e^{\hat\mu t_N}\cos (\hat \omega t_N + \phi)]^T
\end{equation}
Here, $C_1$ is an $\mathcal{O}(1)$ constant, $\phi$ is a phase-shift introduced so that $\mathscr{Y}_m$ remains real-valued, $\hat\lambda^\pm = \hat \mu \pm i \, \hat\omega$ is a known high-order eigenvalue satisfying $|\hat\mu |> |\mu_i|$ and $|\hat \omega| > |\omega_i|$ for each $ i = 1, \cdots, M$, and $\mathscr{Y}_m$ is an error term from the Monte Carlo simulation.

Let $\delta_i$ be the volume of the bin around $\mathbf{x}_i$.\footnote{We assume uniform volume about each sample point $\mbx_i$: $\delta\equiv \delta_i$ for each $\mbx_i\in \mathcal{X}$.}
Assume $K$ samples are drawn from $X_{t_j}$ for each $j = 1, \cdots, N$. Then the Monte Carlo estimate of $\rho_{f,i}(t_j)$ is 
\begin{equation}
    \tilde{\rho}_{f,i}(t_j) = \frac{1}{K\delta_i}B(K, \rho_{f,i}(t_j)\delta_i)
\end{equation}
where $B(K, \rho_{f,i}(t_j)\delta_i)$ is a binomial random variable. 
If $K$ is sufficiently large, $\tilde{\rho}_{f,i}(t_j)$ may be approximated by a normal random variable with mean $\rho_{f,i}(t_j)$ and variance $\rho_{f,i}(t_j)(1 - \rho_{f,i}(t_j))/(K \delta_i)$. Therefore, the $j$-th entry of $\mathscr{Y}_m$ is approximated by
\begin{equation}
    (\mathscr{Y}_m)_j \approx \frac{\rho_{f,i}(t_j)(1 - \rho_{f,i}(t_j))}{K \delta}\varepsilon_j
\end{equation}
where each $\varepsilon_j$ is a standard normal random variable. Since each $\rho_{f,i}$ are $\mathcal{O}(1)$ quantities, we can write $\mathscr{Y}_m = D {\mathbf{\varepsilon}}$ for an $\mathcal{O}(1)$ diagonal matrix $D$ whose $j$-th entry is $\frac{\rho_{f,i}(t_j)(1 - \rho_{f,i}(t_j))}{K \delta_i}$. 

\begin{lemma}\label{lemma}
    Assume the real parts of $M$ known eigenvalues satisfy $0 > \mu_1 > \cdots \mu_M$ and $\hat\lambda^\pm = \hat \mu \pm i \, \hat\omega$ is a known high-order eigenvalue satisfying $|\hat\mu |> |\mu_i|$ and $|\hat \omega| > |\omega_i|$ for each $ i = 1, \cdots, M$. 
    Let the time slices $\{t_j\}_{j = 1}^N$ be an equipartition of $[T_s, T_f]$ for $N$ sufficiently large.
    If Monte Carlo error terms $\{\varepsilon_i\}_{i = 1}^N$ are i.i.d. standard normal random variables, then we have
    \begin{equation}
        \| \beta - \beta_0\| \leq \mathcal{O}(|\mu_1 + \hat\mu| ^{-1}e^{(\hat \mu - \mu_1)T_s}) + \mathcal{O}(N^{-1/2} e^{-\mu_1 T_s})
    \end{equation}
\end{lemma}
A detailed proof of Lemma \ref{lemma} is provided in the supplement.

Lemma \ref{lemma} shows that it is beneficial to include multiple eigenvalues into the linear regression, even if we use the neural network to solve for only the first eigenfunction. In addition, the choice of $T_s$ and $N$ should be properly balanced. Larger $T_s$ can reduce the error caused by higher modes, but increase the error caused by the Monte Carlo noise. The Monte Carlo error can be reduced by using more time slices. However, one should note that Lemma \ref{lemma} assumes independent error terms at distinct $t_j$. In practice, the dependency among $\{ \varepsilon_j\}_{j = 1}^N$ cannot be ignored if the sampling gap $|t_{j+1} - t_j|$ is too small.

\section{Examples} \label{sec: ml examples}

In this section, we apply our machine learning approach to approximate the SKO eigenfunctions and eigenvalues of several model systems.
More examples and error analysis are provided in the supplement.
When possible, we compare results with solutions generated by finite difference methods (see the supplement for more details).

\subsection{2D Noisy Stuart-Landau Oscillator}

We first study a noisy 2D Stuart-Landau oscillator
\begin{equation}\label{eq: ml : 2d sl}
\begin{split}
    dX &= [-4X(X^2+Y^2-1) + \omega Y]dt + \sqrt{2D}\,dW_1(t)
    \\
    dY &= [-4Y(X^2+Y^2-1) - \omega X]dt + \sqrt{2D}\,dW_2(t)
\end{split}
\end{equation}
where we take $\omega=2$ and $D=0.09473$.
For the ANN approach, we consider a numerical domain $\mathcal{R} = [-2,2]\times [-2,2]$ discretized into $N^2$ boxes, with $N=200$.
The leading SKO eigenvalue (computed via finite differences) is $\lambda_1 = -0.1 + 2i$.
Note that the stationary distribution, i.e., the eigenfunction $P_0$, is known exactly
\begin{equation}
    P_0(x,y) = \frac{1}{N}\exp\left(-\frac{1}{D}(x^2+y^2-1)^2\right), \quad N = \iint_{\mathbb{R}^2} dxdy \; \exp\left(-\frac{1}{D}(x^2+y^2-1)^2\right)
\end{equation}

We first consider the approximation of the forward eigenmode $P_{\lambda_1}$ using our ANN approach.
We run 400,000 trajectories and sample the forward density at $N_t=5000$ time-slices with timestep $d t = 0.001$ and $t_{\text{gap}}=10$.
Each realization has an initial condition drawn uniformly from the box containing the point $(1,0)$.
The reference box for the eigenvalue approximation is $\mathcal{B} = [-2,0]\times [-2,0]$.

\begin{figure}[ht]
    \centering
    \includegraphics[scale=.35]{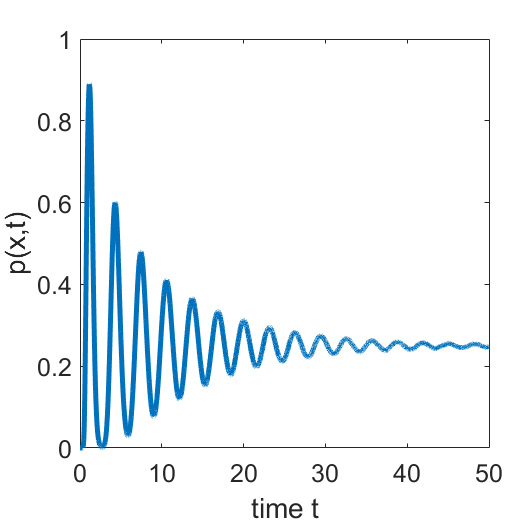}
    \includegraphics[scale=.35]{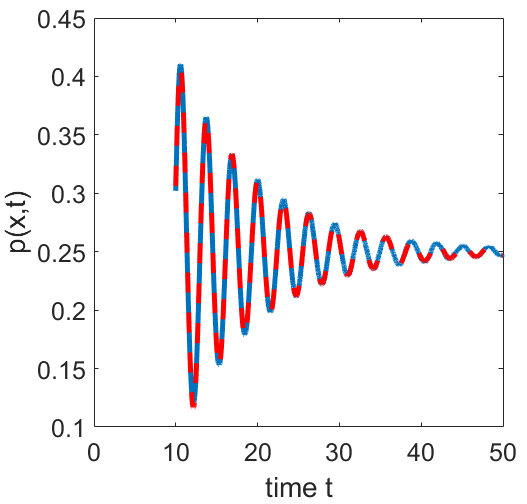}
    
    \includegraphics[scale=.35]{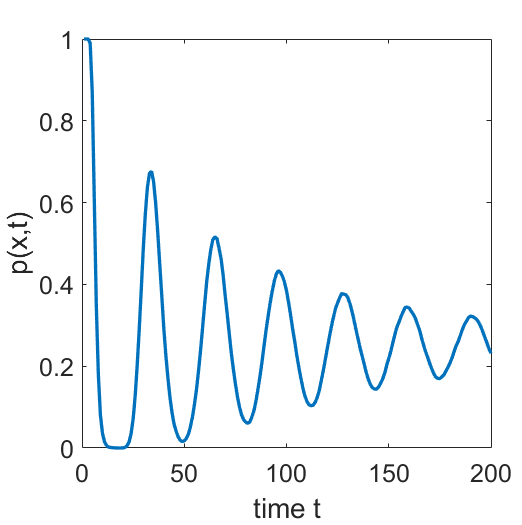}
    \includegraphics[scale=.35]{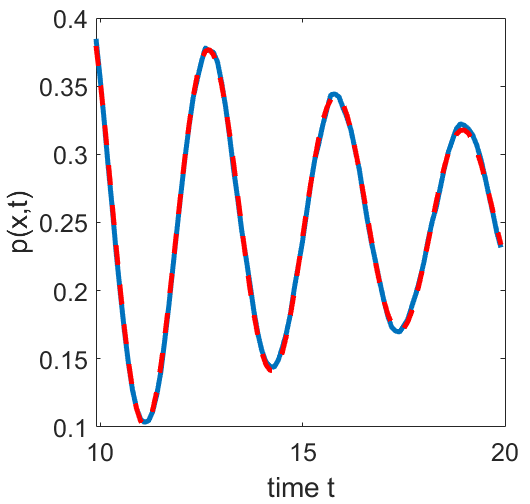}
    \caption{\textbf{Top row:} The decaying probability density associated with the forward eigenmode (left) over $N_t=5000$ time-slices for the 2D Stuart-Landau system \eqref{eq: ml : 2d sl}.
    We fit the decaying density using only time-slices from $t_k$ for $k\in[1000,5000]$ (red, right).
    \textbf{Bottom row:} The decaying probability density associated with the backward eigenmode (left) over $N_t=200$ time-slices for the 2D Stuart-Landau system \eqref{eq: ml : 2d sl}.
    We fit the decaying density using only timeslices from $t_k$ for $k\in[100,200]$ (red, right).}
    \label{fig: ml : 2dsl forward fitting}
\end{figure}

\begin{figure}[ht]
    \centering
    \includegraphics[scale=.35]{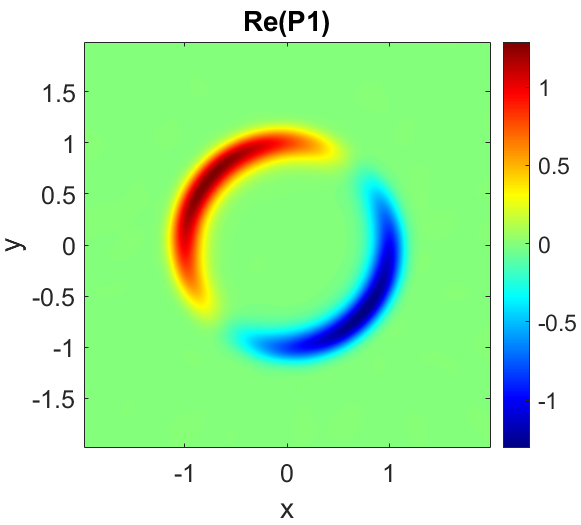}
    \includegraphics[scale=.35]{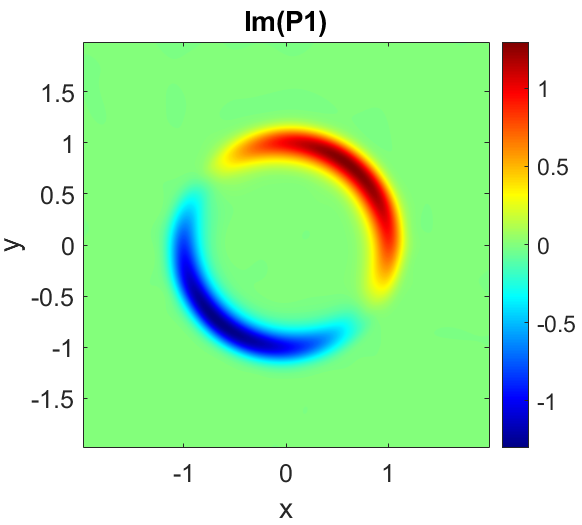}

    \includegraphics[scale=.35]{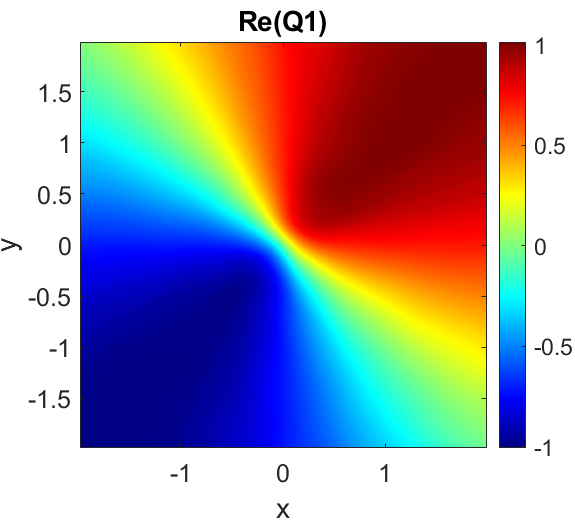}
    \includegraphics[scale=.35]{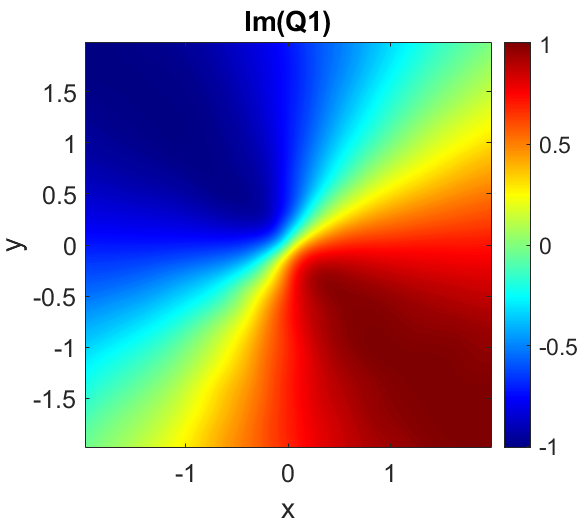}
    \caption{\textbf{Top row:} The ANN approximations of the real (left) and imaginary (right) parts of the forward eigenmode $P_{\lambda_1}$.
    \textbf{Bottom row:} The ANN approximations of the real (left) and imaginary (right) parts of the backward eigenmode $Q^*_{\lambda_1}$.}
    \label{fig: ml : forward 1 mode er}
\end{figure}

In the top row of Figure \ref{fig: ml : 2dsl forward fitting}, we plot the decay of the forward time-dependent density in the reference box.
The decaying density is fit to time-slices $t_k$ with $k\in [1000,5000]$.
Our fitting returns an estimated eigenvalue $\Tilde{\lambda}_1 = -0.0996 + 2.0001 i$.
Both the real and imaginary parts have less than 1\% relative error compared to the FD solution.
We display the output of the ANN in the top row of Figure \ref{fig: ml : forward 1 mode er}.
We used $|\mathcal{Y}|=$500,000 reference points, all 40,000 training points $\mathcal{X}$, and trained the ANN over 30 epochs.

To compute the $Q$-function of the 2D Stuart-Landau model \eqref{eq: ml : 2d sl}, we run 20,000 trajectories with initial conditions drawn uniformly on each of the $N^2=$ 40,000 boxes and sample the backward density at $N_t=200$ time-slices with timestep $d t = 0.005$ and $t_{gap}=20$.
The reference box is $\mathcal{B} = [-2,0]\times [-2,0]$.

In the bottow row of Figure \ref{fig: ml : 2dsl forward fitting}, we plot the decay of (one of) the backward time-dependent densities in the reference box.
The decaying density is fit to time-slices $t_k$ for $k\in [100,200]$.
Our fitting returns an estimate eigenvalue $\Tilde{\lambda}_1 = -0.0992 + 1.9997 i$.
Once again, both the real and imaginary parts have less than 1\% relative error compared to the FD solution.
We also show the output of the ANN in the bottom row of Figure \ref{fig: ml : forward 1 mode er}.
We used $|\mathcal{Y}|=$500,000 reference points, all 40,000 training points $\mathcal{X}$, and trained the ANN over 30 epochs.
Details on the error associated with the approximation of the forward and backward eigenfunctions may be found in the supplement.

\subsection{4D Noisy Coupled Morris-Lecar Neurons}

Here, we study a system of two noisy Morris-Lecar neurons coupled by a gap junction
 \begin{equation}\label{eq: ml : 4D ML}
    \begin{split}
        dV_1 &= \frac{\iota}{C}[I - g_L(V_1-v_L)-g_KN_1(V_1-v_K)-g_{CA}m_\infty(V_1)(V_1-v_{CA})
        \\
        &\quad + \kappa(V_2-V_1)]dt + \sqrt{2D_{v_1}}\,dW_1(t)
        \\
        dN_1 &= \iota\left[\alpha(V_1)(1 - N_1) - \beta(V_1)N_1\right]dt + \epsilon_1\sqrt{\alpha(V_1)(1 - N_1) + \beta(V_1)N_1}\,dW_2(t)
        \\
        dV_2 &= \frac{\iota}{C}[I - g_L(V_2-v_L)-g_KN_2(V_2-v_K)-g_{CA}m_\infty(V_2)(V_2-v_{CA})
        \\
        &\quad + \kappa(V_1-V_2)]dt + \sqrt{2D_{v_2}}\,dW_3(t)
        \\
        dN_2 &= \iota\left[\alpha(V_2)(1 - N_2) - \beta(V_2)N_2\right]dt + \epsilon_2\sqrt{\alpha(V_2)(1 - N_2) + \beta(V_2)N_2}\,dW_4(t)
    \end{split}
\end{equation}
where functions and parameter values for the noiseless system are provided in the supplement, and $\iota=20$.
The noise parameters are $D_{v_1}=15$, $D_{v_2}=50$, $\epsilon_1=0.3$, and $\epsilon_2 = 1$.
We increase the noise of the second oscillator so that its isolated eigenvalues have more negative real parts, i.e., so that the dominant SKO eigenvalues of \eqref{eq: ml : 4D ML} correspond only to the first oscillator.
We set $\kappa=0$ so that the oscillators are uncoupled, i.e., so that we can project a 4D solution onto a 2D plane and compare with known solutions generated by finite differences.
For reference, the first two eigenvalues (computed with finite-differences) are $\lambda_1 = -0.0879 + 1.9518i$ and $\lambda_2 = -0.2521 + 3.9602i$.

\begin{figure}[ht]
    \centering
    \includegraphics[scale=.35]{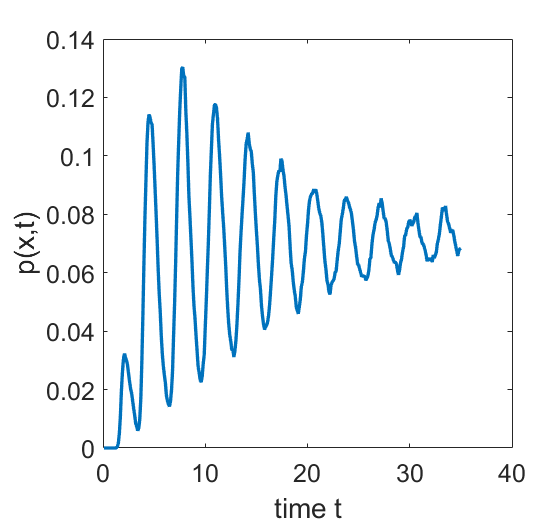}
    \includegraphics[scale=.35]{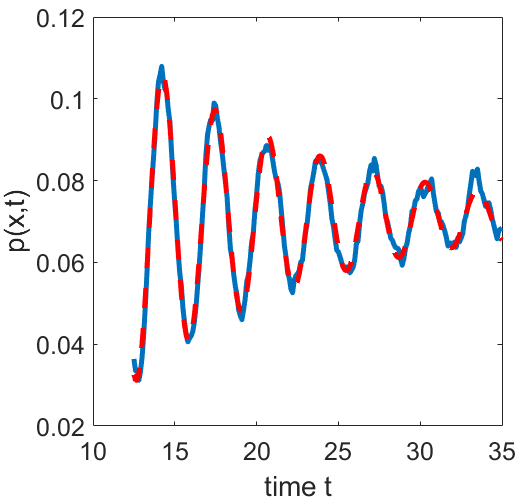}

    \includegraphics[scale=.35]{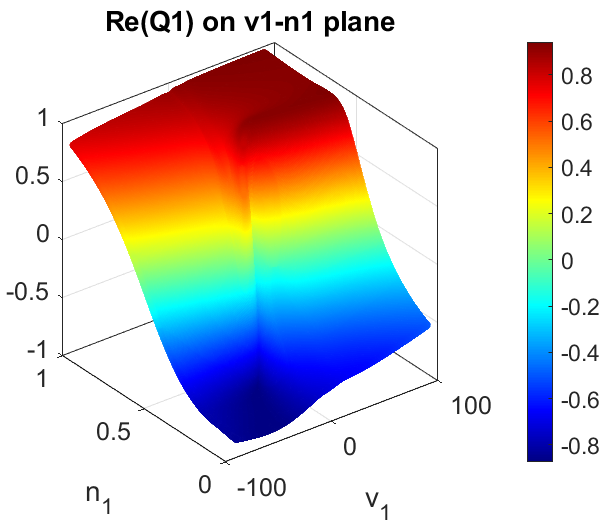}
    \includegraphics[scale=.35]{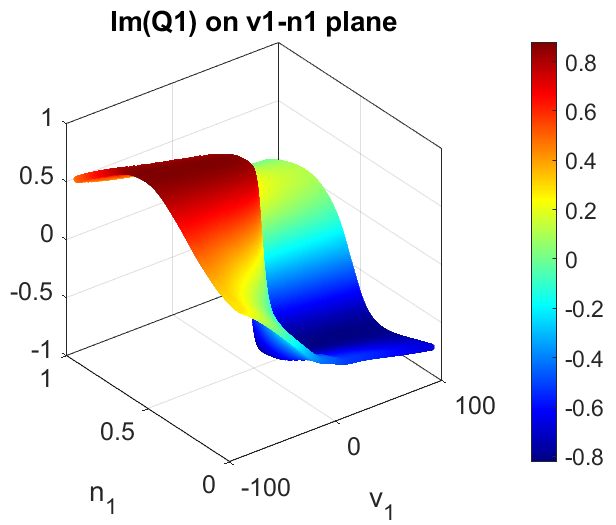}
    \caption{\textbf{Top row:} The decaying probability density (left) over $N_t=350$ time-slices for the 4D Morris-Lecar system \eqref{eq: ml : 4D ML}.
    We fit the decaying density using only timeslices from $t_k$ for $k\in[125,350]$ (red, right).
    \textbf{Bottow row:} ANN approximations of the real (left) and imaginary (right) parts of the $Q$-function, projected onto the $v_1$,$n_1$ plane with $v_2=n_2=0$.}
    \label{fig: ml : 4dml backward fitting}
\end{figure}

Here, we focus only on the computation of the $Q$-function.
To generate densities, we discretize the domain $\mathcal{R} = [-90,90]\times [0,1]\times [-90,90]\times [0,1]$ into $N^4$ boxes, with $N=100$.
We store only 10,000 of the boxes in memory, chosen randomly and uniformly throughout $\mathcal{R}$.
We take 30,000 initial conditions in each box, and keep track of the corresponding densities at $N_t=350$ time-slices.
We set $t_{\text{gap}}=20$ and $dt=0.005$.
In the top row of Figure \ref{fig: ml : 4dml backward fitting}, we depict the decaying probability density corresponding to one of the boxes.
Our fitting procedure estimates that $\Tilde{\lambda}_1 = -0.0871 + 1.9573 i$ taking time-slices $t_k$ with $k\in[125,350]$, which is in excellent agreement with the true value.

To generate approximations of the $Q$-function, we solve the least squares problem \eqref{eq: ml : backward equation} and include the first two non-trivial eigenmodes in the expansion.
We solve only with time-slices $t_k$ for $k\in [125,350]$.
Once this data is generated, we train an ANN on the 10,000 collocation points $\mathcal{X}$,
and $|\mathcal{Y}|=$1,000,000 reference points.
The ANN was trained over 200 epochs, and iterations were divided into 128 batches, which were optimized separately to improve speed.
In this case, the hidden layers of the ANN consist of only 12 nodes.
We rescale voltage $\overline{v_j} = V_j/100$ for $j=1,2$ to be on the same order of magnitude as the gating variables $N_j$ when training the ANN.
The results are shown in the bottom row of Figure \ref{fig: ml : 4dml backward fitting}.

\subsection{Lorenz System with Additive Noise}

As our final example, we consider a 3D Lorenz system in the chaotic regime that is subject to additive noise
\begin{equation}\label{eq: ml : lorenz}
    \begin{split}
        dX &= [\sigma(Y-X)]dt + \sqrt{2D}\,dW_1(t)
        \\
        dY &= [X(\rho -Z) - Y]dt + \sqrt{2D}\,dW_2(t)
        \\
        dZ &= [XY - \beta Z]dt + \sqrt{2D}\,dW_3(t)
    \end{split}
\end{equation}
where $\sigma=10$, $\beta=8/3$, $\rho=28$, and $D=5$.

We are not aware of any solution methods for the SKO eigenfunctions of \eqref{eq: ml : lorenz}, apart from our machine learning approach.\footnote{See these articles for a Koopman perspective on the deterministic Lorenz system \cite{brunton2021modern, jin2023invertible,korda2020data,otto2021koopman, schwabedal2012optimal}.}
The supplement contains information regarding the approximation of the stationary distribution of \eqref{eq: ml : lorenz} when $D\neq0$.

\begin{figure}[ht]
    \centering
    \includegraphics[scale=.35]{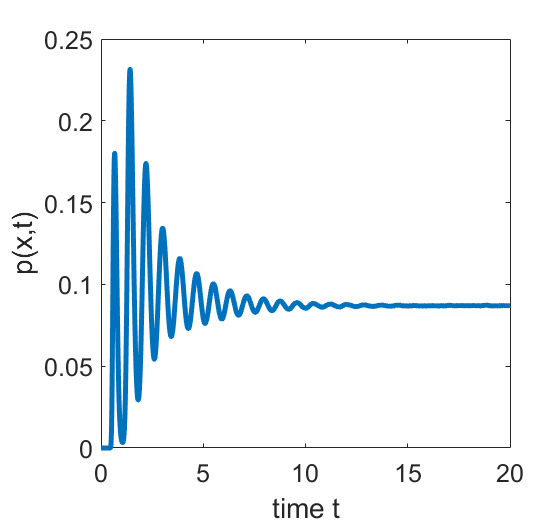}
    \includegraphics[scale=.35]{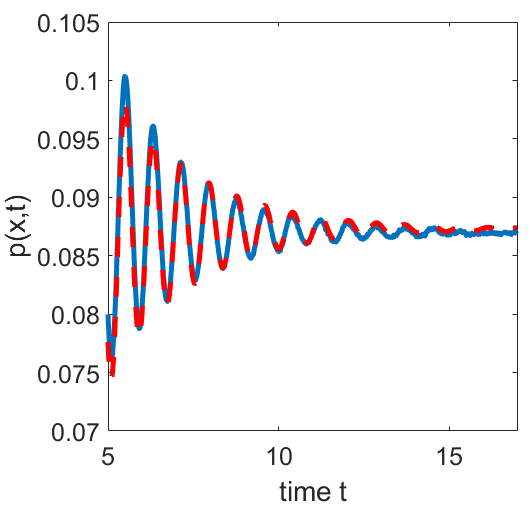}

    \includegraphics[scale=.35]{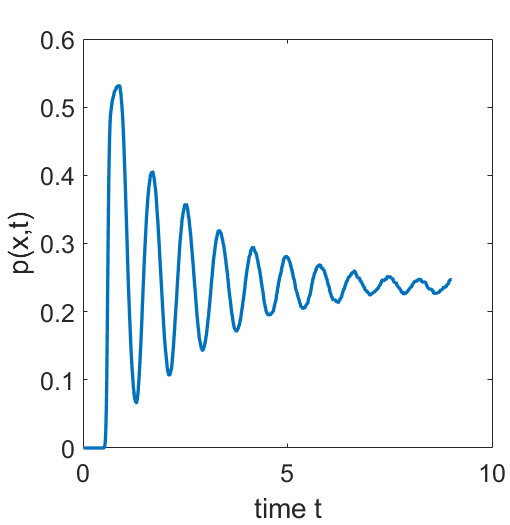}
    \includegraphics[scale=.35]{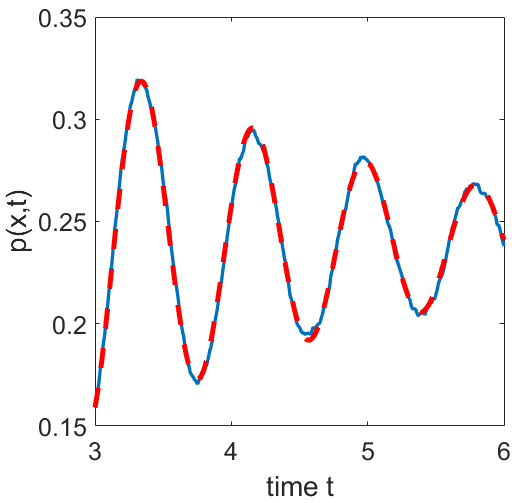}
    \caption{\textbf{Tow row:} The decaying probability density (left) over $N_t=2000$ time-slices for the noisy Lorenz system \eqref{eq: ml : lorenz} with $D=5$.
    We fit the decaying density using only timeslices from $t_k$ for $k\in[500,1700]$ (red, right).
    \textbf{Bottow row:} The decaying probability density (left) over $N_t=900$ time-slices for the noisy Lorenz system \eqref{eq: ml : lorenz}.
    We fit the decaying density using only timeslices from $t_k$ for $k\in[300,600]$ (red, right).}
    \label{fig: ml : lorenz forward fitting}
\end{figure}

\begin{figure}[ht]
    \centering
    \includegraphics[scale=.33]{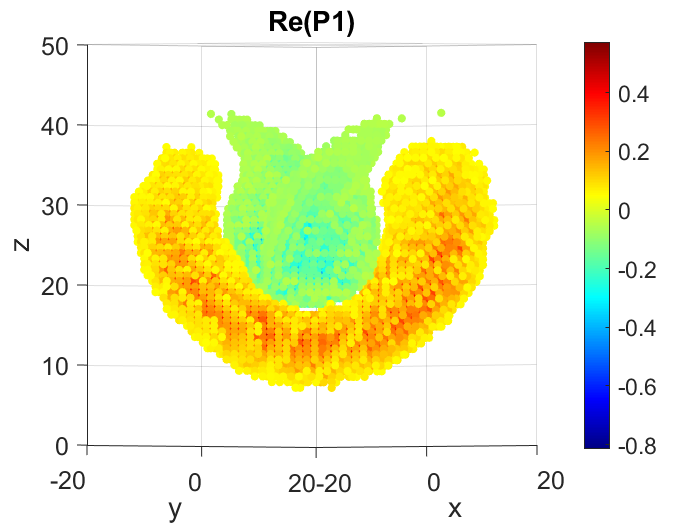}
    \includegraphics[scale=.33]{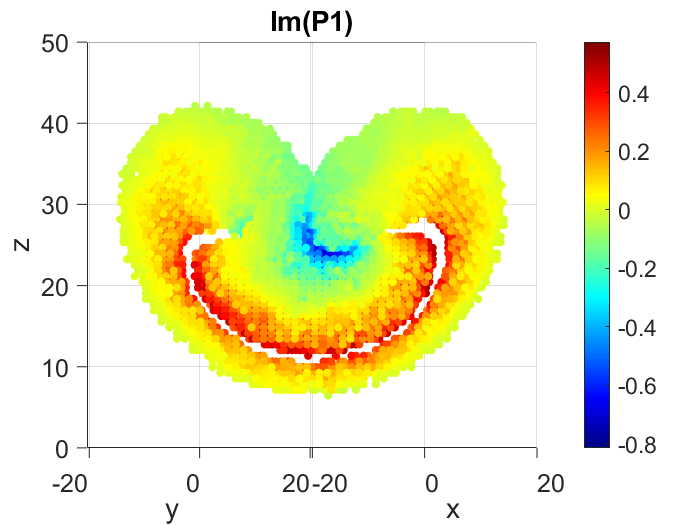}

    \includegraphics[scale=.33]{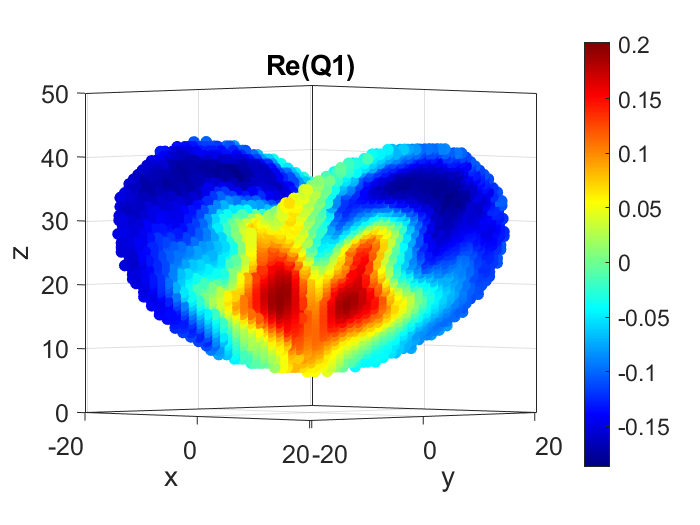}
    \includegraphics[scale=.33]{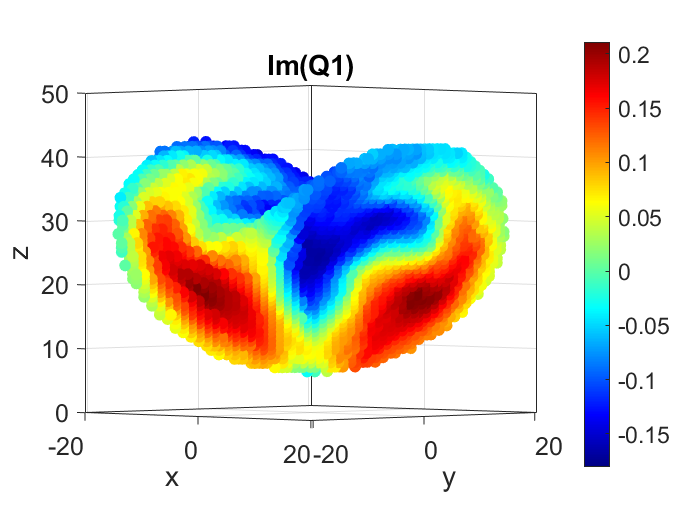}

    \includegraphics[scale=.33]{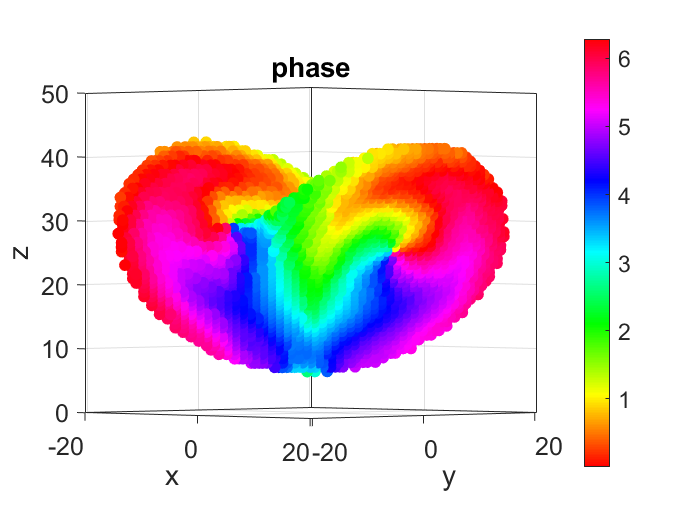}
    \caption{.
    \textbf{Top row:} ANN approximation of the slowest decaying forward eigenmode of the noisy Lorenz system.
    \textbf{Middle row:} Ann approximation of the $Q$-function for the noisy Lorenz system.
    In both cases, any values falling below a threshold of $0.05$ were disregarded to elucidate the structure of the eigenfunctions.
    \textbf{Bottow row:} The ``stochastic asymptotic phase'' \cite{thomas2014asymptotic}, i.e., the complex argument of the ANN approximation of $Q_{\lambda_1}^*$.}
    \label{fig: ml : lorenz forward}
\end{figure}

We now consider the approximation of the forward eigenmode, $P_{\lambda_1}$, using our ANN approach.
We run 10,000,000 trajectories and sample the forward density in 20,000 boxes at $N_t=2000$ time-slices with timestep $d t = 0.001$ and $t_{gap}=10$.
Each realization is initialized with initial condition drawn uniformly from the box that contains the point $(10,10,10)$.
The reference box is $\mathcal{B} = [-30,0]\times [-30,0] \times [-30, 20]$.
In the top row of Figure \ref{fig: ml : lorenz forward fitting}, we show the decay of the forward time-dependent density in the reference box.
The decaying density is fit only for time-slices $t_k$ for $k\in [500,1700]$.
Our fitting returns an estimated eigenvalue $\Tilde{\lambda}_1 = -0.4107 + 7.6684
 i$.
We find the least squares approximation for the forward eigenfunction via \eqref{eq: ml : forward equation}, including only one eigenmode in the computation.
We chose the same time interval as for the eigenvalue approximation, $k\in [500,1700]$.
We also train an ANN, and show its output in the top row of Figure \ref{fig: ml : lorenz forward}.
The ANN was trained over 30 epochs.
We used $|\mathcal{Y}|=$500,000 reference points, and all 20,000 training points $\mathcal{X}$.

Finally, we consider the $Q$-function of the noisy Lorenz system \eqref{eq: ml : lorenz}.
We store 20,000 boxes in the domain $\mathcal{R}$.
We take 20,000 initial conditions in each box, and keep track of the corresponding densities at $N_t=900$ time-slices.
We set $t_{\text{gap}}=100$ and $dt=0.0001$.
In the bottom row of Figure \ref{fig: ml : lorenz forward fitting}, we depict the decaying probability density corresponding to one of the boxes.
Our fitting procedure estimates that $\tilde{\lambda}_1\approx -0.4145 + 7.6709i$ taking time-slices $t_k$ with $k\in [300,600]$.\footnote{This estimate of $\lambda_1$ agrees to two decimal places with the approximation of the forward eigenmodes.}
We find the least squares approximation for the $Q$-function via \eqref{eq: ml : backward equation}, including only one eigenmode in the computation.
We chose the same time interval as for the eigenvalue approximation, $k\in [300,600]$.
We train an ANN over 40 epochs, with $|\mathcal{Y}|=$500,000 reference points, and all 20,000 training points $\mathcal{X}$.
Iterations were broken up into 128 batches, which were optimized separately to improve speed.
We interpolate both the least squares and ANN solutions over the ``butterfly attractor", i.e., the high-density region of the stationary distribution.
The results are shown in the middle row of Figure \ref{fig: ml : lorenz forward}.
The bottom row of Figure \ref{fig: ml : lorenz forward} shows the stochastic asymptotic phase $\Psi(\mathbf{x}) = \arg{Q_{\lambda_1}^*(\mathbf{x})}$.

\section{Discussion}

In this paper, we developed a data-driven machine learning framework for computing SKO and Fokker-Planck eigenfunctions. Our approach builds on prior work on data-driven Fokker-Planck solvers \cite{li2018data,zhai2022deep}, which leverage Monte Carlo simulation data to guide either finite difference schemes or neural network training. Unlike the Fokker-Planck equation, whose solution can be directly approximated using Monte Carlo data, the eigenfunctions cannot be sampled in this way. To address this, we extract the forward (resp.~backward) eigenfunctions via nonlinear regression on time series data of the probability density function (resp.~the expectation of test functions). A central challenge in neural network-based eigenfunction solvers is avoiding convergence to the trivial solution. While earlier methods often employ quadrature-based normalization techniques \cite{han2020solving,jin2022physics}, our regression-based approach directly estimates eigenfunction values at numerous collocation points, thereby accelerating convergence and improving solution quality. The method is flexible, requiring little problem-specific adjustments, and is both mesh-free and naturally suited to high-dimensional settings.

The framework presented here also extends existing theory pertaining to stochastic oscillators.
Previous works \cite{Perez2023universal, perez2021isostables,thomas2014asymptotic} numerically studied the eigenfunctions of planar stochastic oscillators using finite differences, yet were not able to fully analyze systems in dimensions $n\geq 3$.
For the first time, the theory developed in this paper allow for direct numerical investigation of higher-dimensional stochastic oscillators and systems of coupled stochastic oscillators, which has previously not been numerically feasible.
A significant contribution of this paper is to show that upon the introduction of stochastic forcing, the 3D Lorenz system becomes a stochastic oscillator (in the sense introduced in \cite{thomas2014asymptotic}).
The probability density exhibits decaying oscillations, and the slowest decaying SKO eigenmode is complex-valued.
The complex argument of this eigenmode gives a well-defined stochastic asymptotic phase for the system.
In contrast, the \textit{deterministic} Lorenz system exhibits chaotic dynamics.
It does not admit a periodic limit-cycle solution, and therefore does not have a well-defined phase reduction.\footnote{However, an ``optimal'' phase description is possible for deterministic chaotic oscillators, including the Lorenz system \cite{schwabedal2012optimal}.}

In the deterministic case, systems with a spiral sink do not admit periodic-limit cycle solutions.
Nevertheless, when perturbed by noise, the system behaves as a stochastic oscillator with a well-defined asymptotic stochastic phase \cite{thomas2019phase}.
Similarly, deterministic systems which admit heteroclinic cycles do not have periodic limit-cycle solutions and do not admit a well-defined phase reduction.
But when perturbed by noise, such a system does exhibit a well-defined $Q$-phase reduction \cite{thomas2014asymptotic}.
As far as we are aware, no previous work has shown that a deterministic chaotic system, when perturbed by noise, exhibits a well-defined $Q$-phase reduction.
Here, for the first time, we demonstrate numerically that the Lorenz system \eqref{eq: ml : lorenz} does have a well-defined $Q$-phase reduction when perturbed by additive noise.
This observation supports the description of the $Q$-function as a universal description of stochastic oscillators in that it is agnostic to the underlying mechanism of oscillation \cite{Perez2023universal}.

Our machine learning tool allows for the study of high-dimensional systems of coupled stochastic oscillators. Several directions can be explored to further improve training quality. First, our sampling technique could be improved by, for example, using a forward-reverse sampler developed in \cite{milstein2004transition}, which can provide better estimate of the probability density function in high-dimensional setting. Second, our approach could benefit from integrating other data-driven techniques related to the Koopman operator. While extended dynamic mode decomposition (EDMD) has been used to compute eigenfunctions of the Fokker-Planck operator \cite{xu2025data, zhao2023data}, its performance depends strongly on the choice of basis functions, or “dictionary”. One promising direction is to incorporate EDMD into neural network training by using the approximated Koopman operator to provide additional regularization.

On the application side, our future work may focus on extending the current framework to compute a family of $Q$-functions and corresponding eigenvalues as a function of coupling strength. Another compelling direction is to adapt the ANN-based framework for systems with partial observations, enabling the computation of eigenfunctions for partially observed oscillators. Computing SKO eigenfunctions is also of significant interest to experimentalists, especially in scenarios where the underlying model is unknown. The framework developed here could be further modified to infer underlying dynamics directly from time-series data, making it applicable to a wide range of real-world problems.

\section*{Acknowledgments}
This work was supported in part by NSF grant DMS-2052109 to PJT and DMS-2108628 to YL. 
This material is also based upon work supported in part by the National Science Foundation under Grant No.~DMS-1929284 while the authors were in residence at the Institute for Computational and Experimental Research in Mathematics in Providence, RI, during the ``Math + Neuroscience: Strengthening the Interplay Between Theory and Mathematics" program.
This work was also supported in part by the Oberlin College Department of Mathematics.

\bibliographystyle{Inputs/siamplain}
\bibliography{Inputs/lib}

\end{document}